\documentstyle{article}
\newcommand{\ncm}{\newcommand}
\ncm{\al}{\alpha}
\ncm{\bt}{\beta}
\ncm{\C}{\bf C}
\ncm{\cc}{C_c(E)}
\ncm{\co}{C_0(E)}
\ncm{\cs}{C^*_r}
\ncm{\cse}{\cs(E)}
\ncm{\dl}{\delta}
\ncm{\Dl}{\Delta}
\ncm{\df}{\partial}
\ncm{\ep}{\epsilon}
\ncm{\gm}{\gamma}
\ncm{\Gm}{\Gamma}
\ncm{\io}{\iota}
\ncm{\K}{{\cal K}}
\ncm{\Ll}{{\cal L}}
\ncm{\lm}{\lambda}
\ncm{\lt}{\langle}
\ncm{\N}{{\cal N}}
\ncm{\ol}{\overline}
\ncm{\om}{\omega}
\ncm{\Om}{\Omega}
\ncm{\op}{\bigoplus}
\ncm{\ot}{\otimes}
\ncm{\pr}{{\em Proof}:\ }
\ncm{\ph}{\varphi}
\ncm{\qef}{\hfill \rule{.3em}{1.6ex}}
\ncm{\ra}{\rightarrow}
\ncm{\rt}{\rangle}
\ncm{\sg}{\sigma}
\ncm{\Sg}{\Sigma}
\ncm{\te}{\theta}
\ncm{\T}{\bf T}
\ncm{\zt}{\zeta}

\begin{document}

\title{Fell bundles over groupoids}
\author{Alex Kumjian}
\date{}
\maketitle

A C*-algebraic bundle (see \cite[\S 11]{F2}) over a locally compact
group may be thought of as a continuous version of a group grading 
in a C*-algebra; one may regard the associated C*-algebra as a fairly 
general sort of crossed product of the fiber algebra over the neutral 
element by the group (in \cite{LPRS} it is shown that the C*-algebra is
endowed with a coaction by the group).  There is a natural extension
of this definition to groupoids (see \cite{Yg}) which when specialized
to trivial groupoids (i.e. topological spaces) yields the more usual 
notion of C*-algebra bundle.  This object is referred to below as a 
Fell bundle.

Closely related notions have appeared in the literature:  in recent
work \cite{Yn}, Yamanouchi studies the analogous notion in the von 
Neumann algebra setting under the name integrable coaction.  Fell 
bundles are reminiscent of the C*-categories discussed in \cite{GLR}.
They are also presaged in \cite[ Def. 5.3]{Re2} Def. 5.3 (the object 
is used to construct a strong Morita equivalence bimodule).

Since each fiber of a saturated Fell bundle may be regarded as a
strong Morita equivalence bimodule, the theory of Fell bundles
provides a natural locus for proving theorems related to strong 
Morita equivalence of the kind which appear in \cite{MRW} and 
\cite{Re2}.

In \S 1  we fix notation and review some well-known facts concerning 
groupoids, Banach bundles, Hilbert modules, and equivalence bimodules.
The notion of Fell bundle is defined in \S 2  and some examples are 
discussed, including the semi-direct product which results from a 
groupoid acting on a C*-algebra bundle fibered over the unit space.  
In \S 3 the associated C*-algebra is constructed in the case that the 
groupoid is r-discrete; the norm is defined by an analog of the left 
regular representation (the resulting C*-algebra may be regarded as 
the reduced C*-algebra associated to the Fell bundle).  Finally, a 
Morita equivalence theorem of the sort discussed above is proved 
(Th. 4.2).  This is used to show that there is an action of a groupoid
on a C*-algebra bundle obtained from a saturated Fell bundle so that 
the C*-algebra associated to the Fell bundle is strong Morita
equivalent to that of the semi-direct product (Cor. 4.5).  See 
\cite[ Th. 8]{Kt}, \cite[Cor. 3.7]{PR}, \cite[Th. 3.1]{Q1}, 
\cite[Cor. 2.7]{Q2}, and \cite[Th. 7.8]{Yn} for related results.

{\bf \S 1  Preliminaries}

{\bf 1.1\ } Given a groupoid $\Gm$, let $\Gm^0$  denote the unit space
and $r, s : \Gm \ra \Gm^0$  denote the range and source maps
(respectively); let  $\Gm^2$  denote the collection of composable pairs; 
write the inverse map $\gm  \mapsto \gm^*$.  We tacitly assume that all 
groupoids under discussion are locally compact and Hausdorff, that the
structure maps are continuous, and that they admit left Haar systems 
(see \cite{Re1}).  Let $\Dl$  denote the transitive equivalence
relation with unit space  $\Dl^0 = \{0,1\}$;  write  
$\Dl = \{ 0, 1, \df , \df^*\}$  where $s(\df) = 0$  and  $r(\df) = 1$.  
A groupoid  $\Gm$  is said to be trivial if  $\Gm = \Gm^0$.

{\bf 1.2\ }  Given a Banach bundle, $ p : E \ra X $, let  $\cc$
denote the collection of compactly supported continuous sections of  
$E$  and $\co$  denote the collection  of continuous sections of $E$ 
vanishing at $\infty$.  Note that $\co$  may be regarded as the
completion of $\cc$ in the supremum norm.  For $x \in X$, let $E_x$  
denote the fiber over $x$ , $p^{-1}(x)$.  We shall tacitly assume
that the total space, $E$, of any Banach bundle under consideration 
is second countable and that the base space, $X$, is locally compact 
and Hausdorff; it follows that the base space is second countable and 
that both the fiber, $E_x$, and the associated Banach space, $\co$,
are separable (see \cite[Prop. 10.10]{F2}).   By a result of Doady and
dal Soglio-H\'{e}rault ( \cite[appendix]{F2}), if  $X$  is locally
compact,  $E$  has enough continuous sections; thus for every 
$e \in E$ there is $f \in \cc$  such that  $f(p(e)) = e$.

{\bf 1.3\ } Let  $A$  be a C*-algebra and $V$ be a right $A$-module;
$V$  is said to be a (right) pre-Hilbert $A$-module if 
it is equipped with an $A$  valued inner product $\lt \cdot, \cdot \rt $  
which satisfies the following conditions:
\begin{itemize}
\item[i]  $\lt u, \lm v + w \rt  = \lm \lt u,  v \rt  + \lt u,  w \rt $
\item[ii] $ \lt u,  va \rt  =  \lt u,  v \rt  a $
\item[iii] $ \lt v,  u \rt  =  \lt u,  v \rt ^*$
\item[iv] $\lt v, v \rt  \ge 0$ and  $\lt v, v \rt  = 0$ only if $v = 0$.
\end{itemize}
for all  $u, v, w \in V$,  $\lm \in \C$, and  $a \in A$.  Say that
$V$  is a (right) Hilbert $A$-module if it is complete in the norm: 
$\|v\| = \| \lt v, v \rt  \|^{1/2}$.
There is an analogous definition for left Hilbert $A$-modules (the
inner product in this case is linear in the first variable).  Given a
right Hilbert $A$-module $V$, there is a left Hilbert $A$-module
$V^*$  with a conjugate linear isometric isomorphism from  $V$  to
$V^*$,  write $v \mapsto v^*$, which is compatible with the module structure 
and inner product in the following way:
$$av^* = (va^*)^*$$
$$\lt u^*, v^*\rt  = \lt u, v \rt $$
for all  $u, v \in V$  and $a \in A $.  The term, Hilbert $A$-module,
will be understood to mean right Hilbert $A$-module.  If the span of
the values of the inner product of a Hilbert $A$-module $V$ is dense 
in  $A$, then $V$ is said to be full 
(see \cite[\S 2]{Ks}, \cite[13.1.1]{B}).  Note that $A$ may be
regarded as a Hilbert $A$-module when 
endowed with the inner product $ \lt a,  b \rt  = a^*b$;
evidently,  A  is full.

{\bf 1.4\ } Given Hilbert $A$-modules $U$  and $V$, let  $\Ll(V, U)$  
denote the collection of bounded adjointable operators from $V$ to $U$
which commute with the right action of $A$.  Put  $\Ll(V) = \Ll(V, V)$; 
endowed with the operator norm $\Ll(V)$ is a C*-algebra.  
For  $u \in U$, $v \in V$, define  $\te_{u,v} \in \Ll(V, U)$  by
$\te_{u,v}(w) = u \lt v, w \rt $ for  $w \in V$.  Note that 
$\te_{u,v}=\te_{v,u}^*$.  
The closure of the span of such operators is denoted $\K(V, U)$;
write $\K(V) = \K(V, V)$.  Note that $\K(V)$ is an (essential) ideal in
$\Ll(V)$; in fact, $\Ll(V)$ may be identified with the multiplier 
algebra of $\K(V)$ (see \cite[Th. 1]{Ks}).

{\bf 1.5\ } Given C*-algebras $A$ and $B$, a $B$-$A$ equivalence
bimodule $V$  is both a full right Hilbert $A$-module and a full left 
Hilbert $B$-module with compatible inner products ( $V^*$  may be
viewed as an $A$-$B$ equivalence bimodule); note that equivalence
bimodules are also known as imprimitivity bimodules 
(see \cite[Def. 6.10]{Ri1}).  If such a bimodule exists the 
two C*-algebras are said to be strong Morita equivalent (see
\cite{Ri2}).   If $V$ is a full Hilbert $A$-module, then $V$ is a 
$\K(V)$-$A$ equivalence bimodule with $\K(V)$ valued inner product:
$\lt u, v \rt  = \te_{u,v}$.  Indeed, every equivalence bimodule is of 
this form.  Moreover, if  $U$  and $V$  are full Hilbert $A$-modules 
then $\K(V, U)$  may be regarded as a $\K(U)$-$\K(V$) equivalence bimodule.

{\bf 1.6\ }  Given C*-algebras $A$, $B$, and $C$  together with a
$B$-$A$ equivalence bimodule $V$  and a $C$-$B$  equivalence bimodule $U$,
one may form the $C$-$A$ equivalence bimodule:
$U \ot_B V$ (see \cite[Th. 5.9]{Ri1}).
Note that if $V$ is a $B$-$A$ equivalence bimodule and $W$ is a $C$-$A$ 
equivalence bimodule then one may identify:
$$W \ot_A V^* = \K(V, W)$$ 
(as $C$-$B$ equivalence bimodules) via the map:
$w \ot v^* \mapsto \te_{u,v}$.
By this identification and the associativity of the tensor product one
obtains $\K(U \ot_B V) \cong \K(U)$:
$$\K(U \ot_B V) = (U \ot_B V) \ot_A (U \ot_B V)^* \cong 
(U \ot_B (V \ot_AV^*)) \ot_B  U^* \cong (U \ot_B B) \ot_B U^* 
\cong U \ot_B U^* = \K(U).$$
Note that this gives an isomorphism of C*-algebras; 
in fact:  $K(U \ot_B V) \cong C \cong \K(U)$.

{\bf 1.7\ }  Given a C*-algebra bundle $A$ over a space $X$, a Banach
bundle $V$  over $X$ is said to be a Hilbert $A$-module bundle if each
fiber $V_x$ is a Hilbert $A_x$-module with continuous module action
and inner product.  Equipped with the natural inner product $C_0(V)$
is a Hilbert $C_0(A)$-module and it is full if and only 
if  $V_x$  is full for every  $x \in X$ (in this case $V$ is said to 
be full).  Associated to $V$ one obtains a C*-algebra bundle $\K(V)$
where $\K(V)_x = \K(V_x)$.  One has, $C_0(\K(V)) \cong \K(C_0(V))$. 

{\bf \S 2  Fell bundles}

We define below the natural analog of Fell's C*-algebraic bundles 
(cf.\ \cite{F2}, \cite{Yg}) for groupoids.  This notion is a
generalization of both C*-algebraic bundles (over groups) and 
C*-algebra bundles (over 
spaces).

{\bf 2.1\ } Let $\Gm$ be a groupoid and $p : E \ra \Gm$  a 
Banach bundle; set
 $$E^2 = \{ (e_1, e_2) \in E \times E : (p(e_1), p(e_2))  \in \Gm^2 \} .$$
{\bf Def:\ } A {\em multiplication} on $E$ is a continuous map,  
$E^2 \ra E$, write $(e_1, e_2) \mapsto e_1e_2$, which satisfies:
\begin{itemize}
\item[i]  $p(e_1e_2) = p(e_1)p(e_2)$  for all  $(e_1, e_2) \in E^2$
\item[ii]  the induced map,  
           $E_{\gm_1} \times  E_{\gm_2} \ra  E_{\gm_1\gm_2}$, 
         is bilinear for all  $(\gm_1, \gm_2) \in \Gm^2$
\item[iii] $(e_1e_2)e_3 = e_1(e_2e_3)$  whenever the multiplication is defined
\item[iv]  $\|e_1e_2\| \le \|e_1\|\,\|e_2\|$ for all $(e_1, e_2) \in E^2 $.
\end{itemize}
An {\em involution} on $E$ is a continuous map, $E \ra E$, write  
$e \mapsto e^*$, which satisfies:
\begin{itemize}
\item[v]  $p(e^*) = p(e)^*$  for all $e \in E$
\item[vi]  the induced map, $E_{\gm} \ra E_{\gm^*}$, is conjugate
          linear for all $\gm \in \Gm$
\item[vii] $ e^{**} = e$ for all $ e \in E$.
\end{itemize}
Finally, the bundle $E$ together with the structure maps is said to be
a Fell bundle if in addition the following conditions hold:
\begin{itemize}
\item[viii] $ (e_1e_2)^* = e_2^*e_1^* $ for all $(e_1, e_2) \in E^2$
\item[ix] $ \|e^*e\| = \|e\|^2$  for all  $e \in E$
\item[x] $e^*e \ge 0$  for all $e \in E$.
\end{itemize}
Note that if $x \in \Gm^0$ then $E_x$ is a C*-algebra (with norm, 
multiplication, and involution induced from the bundle); 
if  $e \in E_{\gm}$  then $e^*e \in E_{s(\gm)}$, hence it makes 
sense to require that $e^*e$  be positive.  Yamagami refers to such 
a bundle as a C*-algebra over a groupoid (see \cite{Yg}).

{\bf 2.2\ } $E$  is said to be nondegenerate if $E_{\gm} \ne {0}$ 
 for all $\gm \in \Gm$; note that  $E_{\gm}$ is a right Hilbert  
$E_{s(\gm)}$-module with inner product  $\lt e_1, e_2\rt  = e_1^*e_2$ 
and a left Hilbert $E_{r(\gm)}$-module with inner product  
$\lt e_1, e_2\rt  = e_1e_2^*$.  Note also that  $E_{\gm}^* \cong E_{\gm^*}$.

{\bf 2.3\ } Given a Fell bundle $E$ over $\Gm$, let $E^0$ denote the 
restriction $E|_{\Gm^0}$;  clearly, $E^0$ is a C*-algebra 
bundle and  $C_0(E^0)$  is a C*-algebra (with pointwise operations).

{\bf 2.4\ } The Fell bundle $E$ is said to be saturated if  
 $E_{\gm_1}\cdot E_{\gm_2}$ is total in  $ E_{\gm_1\gm_2}$ 
for all $(\gm_1, \gm_2) \in \Gm^2$; note that if $E$  is saturated 
then  $E_{\gm}$ may be regarded as an $E_{r(\gm)}$-$E_{s(\gm)}$  
equivalence bimodule (with inner products as above).  
For  $(\gm_1, \gm_2) \in \Gm^2$, one has  
 $E_{\gm_1} \ot_{E_x}  E_{\gm_2} \cong E_{\gm_1\gm_2}$ 
where $ x = s(\gm_1) = r(\gm_2)$  
(via the map  $e_1 \ot e_2 \mapsto e_1e_2$).

{\bf 2.5\ } Examples:
\begin{itemize}
\item[i]
Let $E$  be a C*-algebra bundle over a space $X$.  If we regard $X$
as a trivial groupoid (so $X = X^0$), then $E$ is seen to satisfy 
the above definition.
\item[ii]
Let $E$ be a C*-algebraic bundle over a locally compact group $G$ 
(in the sense of Fell).  If we regard  $G$  as a groupoid 
(so  $G^0 = {1_G}$ ) then  $E$  is easily seen to satisfy the 
above definition (it is essentially the same definition).
\item[iii]  
Let $A$ and $B$ be  C*-algebras and let $C$ be an $B$-$A$ 
equivalence bimodule.  Form a Fell bundle $E$  
over the groupoid $\Dl$  as follows:  set 
$E_0 = A$,  $E_1 = B$, $E_{\df} = C$,$E_{\df^*} = C^*$;
since $\Dl$ is discrete the Banach bundle structure is trivial.  
One defines multiplication and involution in the obvious way and 
checks that the above conditions are satisfied.  Note
that  E  is saturated.
\item[iv]
Let $\Sg$  be a proper $\T$-groupoid over $\Gm$ (see
\cite[Def. 2.2]{Ku}).   Form the associated line bundle:
$E = \Sg *_{\T}\C = (\Sg \times \C) / \T $ (where  
$t(\sg, z) = (ts, t^{-1}z)$ ).
One defines multiplication and involution as follows:
$$(\sg_1, z_1)(\sg_2, z_2) = (\sg_1\sg_2, z_1z_2),$$ 
$$(\sg, z)^* = (\sg^*,\ol z );$$
one must also check that both are well-defined (cf.\ \cite[\S 12]{F2}).
\item[v]
An action of $\Gm$ on a C*-algebra bundle,  $q : A \ra \Gm^0$, is a 
continuous map (see \cite{Re2}, \cite{M}): $\al : \Gm * A \ra A$  
(where $\Gm * A  = \{ (\gm, a) \in \Gm \times A : s(\gm) = q(a) \}$ ),
write $\al(\gm, a) = \al_{\gm}(a)$, which satisfies the following conditions:
\begin{itemize}
\item[a] $ q(\al_{\gm}(a)) = r(\gm)$  for all $\gm \in \Gm$  and  
         $a \in A_{s(g})$
\item[b] $\al_{\gm} : A_{s(g)} \ra  A_{r(g)}$ is a *-isomorphism for
         all  $\gm \in \Gm$ 
\item[c]  $\al_x(a) = a $ for all $x \in \Gm^0$  and $a \in A_{s(g)}$
\item[d]  $\al_{\gm_1\gm_2}(a) = \al_{\gm_1}(\al_{\gm_2}(a))$  
         for all $(\gm_1,\gm_2) \in \Gm_2$  and  $a \in A_{s(\gm_2}$) .
\end{itemize}
Form the semi-direct product (cf.\ \cite[\S 12]{F2})
$\Gm \times_{\al} A = \Gm * A $,
with multiplication and involution given by the formulas:
$$(\gm_1, a_1)(\gm_2, a_2) = (\gm_1\gm_2, \al_{\gm_2^*}(a_1)a_2),$$
$$(\gm, a)^* = (\gm^*, \al_{\gm}(a^*)).$$
Note that as a Banach bundle, $\Gm \times_{\al} A$ is the pull-back of 
$A$  by $s$.  It is routine to check that with this norm and the above
operations, the semi-direct product, $\Gm \times_{\al} A$,
is a Fell bundle over $\Gm$.
\item[vi]
Let $H$  be a Hilbert bundle over a locally compact space $X$ with an 
inner product  $\lt \cdot, \cdot \rt $  which is conjugate linear in the 
first variable.  Let $\Gm$  denote the transitive groupoid  
$X \times X$ (with obvious structure maps).  There is a natural Fell 
bundle over $\Gm$  associated to $H$.  Set 
$E_{(x,y)} = \K(H_y, H_x)$.  The topology of $E$ is prescribed by
giving a linear space of norm continuous sections which is dense in 
every fiber (see \cite[10.4]{F2});  for each pair of continuous
sections, $\xi, \eta$, of $H$  define a continuous section 
$\te_{\xi,\eta}$  of $E$ by the formula:
$$\te_{\xi,\eta}(x,y)\zt = \xi\lt\eta(y), \zt\rt $$  
for all $(x, y) \in \Gm , \zt \in H_y$.
The span of such sections determine a bundle topology for $E$.  
Multiplication is given by composition and 
involution by the usual adjoint.
\end{itemize}

{\bf 2.6 Remark:\ } Let $j : \Om \ra \Gm$  be a continuous groupoid morphism 
and  $p : E \ra \Gm$  be a Fell bundle.  The pull-back of $E$ by $j$,
$j^*(E) = \Om * E = \{(\om, e) \in \Om \times E : j(\om) = p(e)\}$,
may be regarded as a Fell bundle over $\Om$ in the obvious way.

{\bf \S 3  Construction of the associated C*-algebra}

{\bf 3.1\ } Assume that $\Gm$  is an r-discrete groupoid 
(see \cite{Re1}).  Let  $p : E \ra \Gm$   be a Fell bundle; 
we construct the analog of the reduced C*-algebra $\cse$  as 
in \cite[\S 2]{Ku} (in \cite{Yg} the full C*-algebra is constructed).  
Given  $f, g \in \cc$  define multiplication and involution by 
means of the formulas:
$$fg(\gm) = \sum_{\gm=\al\bt}  f(\al)g(\bt) ,$$
$$f^*(\gm) = f(\gm^*)^*.$$
With these operations $\cc$ forms a *-algebra.  Let  
$P : \cc \ra C_c(E^0)$  be the restriction map.  Define a  
$C_c(E^0)$-valued inner product on $\cc$  by:
$\lt f, g \rt  = P(f^*g)$.

{\bf 3.2  Prop:\ }  With this inner product, $\cc$ is a pre-Hilbert  
$C_0(E^0)$-module.\\[.5ex]
\pr  We verify that  $\lt f, f \rt $ is positive as an element of the 
C*-algebra $C_0(E^0)$  for every $f \in \cc$:
$$\lt f, f \rt (x) = f^*f(x) = \sum_{x = \al\bt}  f^*(\al)f(\bt) 
= \sum_{x = s(\gm)} f(\gm)^*f(\gm) \ge 0 $$
for all  $x \in \Gm^0$, and $f \in \cc$; the same computation shows 
that if  $\lt f, f \rt $ = 0  then  $f = 0$.  The remaining 
properties are left for the reader to verify. \qef

For $f \in \cc$, put $\|f\|_2 = \|\lt f, f \rt \|^{1/2}$ and denote the 
completion of $\cc$ in this norm by $L^2(E)$ (which is now a Hilbert  
$C_0(E^0)$-module).

{\bf 3.3\ } We show that $L^2(E)$ is the Hilbert module associated to 
a bundle of Hilbert modules over $\Gm^0$ (see 1.7).  For each  
$x \in \Gm^0$  consider the Hilbert $E_x$-module:
$V_x = \op_{x = s(\gm)} E_{\gm}$,
with inner product:
$$\left\lt  \sum_{x = s(\gm)} c_\gm, \sum_{x = s(\gm)} d_\gm  \right\rt  = 
 \sum_{x = s(\gm)} {c_\gm}^*d_\gm.$$
One may obtain a bundle topology on the union of the fibers by using 
elements of $\cc$ to provide continuous sections in the obvious way 
(given $f \in \cc$  one obtains the section 
$x \mapsto \sum_{x = s(\gm)} f(\gm) \in V_x$); 
let $V$ denote the bundle obtained in this way.  Then $V$ is a 
Hilbert $E^0$-module bundle; one has $C_0(V) \cong L^2(E)$ as Hilbert 
$C_0(E^0)$-modules.  Since $E^0 \subset V$, as Hilbert $E^0$-module 
bundles, $V$ is full.  If $E$ is saturated then for all  
$\bt \in \Gm$  one has:
$$V_{r(\bt)} \ot_ {E_{r(\bt)}} E_\bt \cong V_{s(\bt)} $$ 
via the identification:
$$\left( \sum_{s(\gm)=r(\bt)} c_\gm\right) \ot e \mapsto
 \sum_{s(\gm)=r(\bt)} c_\gm e $$
where $e \in E_\bt$, $c_\gm \in E_\gm$, and 
$$\sum_{s(\gm)=r(\bt)} c_\gm \in \op_{s(\gm)=r(\bt)} E_{\gm} = V_{r(\bt)}.$$
Moreover, one has:
$${V_x}^* = \op_{x = r(\gm)} E_{\gm},$$
(via the involution map) and:
$$E_\bt \ot_{E_{s(\bt)}} {V_{s(\bt)}}^* \cong  {V_{r(\bt)}}^*.$$ 

{\bf 3.4\ } Left multiplication by an element in $\cc$  is a bounded 
operator with respect to the norm, $\|\cdot\|_2$, 
and hence extends to the completion.  One checks :
$$ \lt fg, h \rt  = \lt g, f^*h \rt  $$  
for all $f, g, h \in \cc$; hence, left multiplication by an element 
in $\cc$ is adjointable and one obtains a *-monomorphism:
$$\cc \ra \Ll(L^2(E)).$$ 
Let $\cse$  denote the completion of $\cc$  with respect to the
operator norm; since $\cse$ is a closed *-subalgebra of $\Ll(L^2(E))$,
it is a C*-algebra.  Moreover, for each $x \in \Gm^0$ one has 
a representation:
$$\pi_x : \cse  \ra \Ll(V_x),$$
so that for each $a \in \cse$, $\|a\| = \sup \|\pi_x(a)\|$.
Note:  Every bounded continous complex-valued function $g$ on 
$\Gm^0$ may be identified with an element of 
the multiplier algebra,  $M(\cse)$, as follows:  for  
$f \in \cc$  put $gf(\gm) = g(r(\gm))f(\gm)$; one checks that this 
defines an element of $\Ll(L^2(E))$ which centralizes $\cse$.

{\bf 3.5 Examples:\ }
\begin{itemize}
\item[i]
If $E$ is a Fell bundle over a trivial groupoid $X$ (so $X = X^0$ ), 
then $E$ is a C*-algebra bundle and $\cse = \co$.
\item[ii]
Refer to example 2.5iii above; if $C$ is a $B$-$A$ equivalence
bimodule and $E$ is the associated Fell bundle over $\Dl$, then $\cse$  
is the linking algebra associated to $C$.
(cf.\ \cite[Th. 1.1]{BGR}).
\item[iii]
Refer to example 2.5iv and assume that $\Gm$ is a principal r-discrete
groupoid.  Then $\cse$ has a diagonal subalgebra (see \cite[\S 2]{Ku}) 
isomorphic to $C_0(G^0)$; the twist invariant for the diagonal pair,
$(\cse,C_0(G^0))$, is the inverse of $[\Sg]$. 
\item[iv]
Let  $\Gm$ be a transitive equivalence relation on a countable set 
(with the discrete topology) and $E$ be a saturated Fell bundle over
$\Gm$.   Let $V$ be the Hilbert $E^0$-module bundle over $G^0$
described above (3.3); for every $x \in \Gm^0$, $V_x$  is full so
$V_x$ is a $\K(V_x)$-$E_x$ equivalence bimodule (see 1.5).  Now, since
$E$ is saturated, it follows that for any $\gm \in \Gm$ with  
$r(\gm)= x$ and  $s(\gm) = y$, one has 
$V_{x} \ot_ {E_{x}} E_\gm \cong V_{y} $ (see 3.3).  
This induces a *-isomorphism  $\K(V_y) \cong \K(V_x)$ (see 1.6).  
Moreover, for each $x \in G^0$, $\pi$  induces a *-isomorphism:
$$\cse \cong \K(V_x).$$
\end{itemize}

{\bf  3.6 Prop:\ }  The restriction map, $P: \cc \ra C_c(E^0)$, 
extends to a conditional expectation   $P: \cse \ra C_0(E^0)$. \\[.5ex] 
\pr  This follows from Tomiyama's characterization of a conditional 
expection as a projection of norm one onto a subalgebra (see \cite{T}).  
One checks that for $f \in \cc$:
$$\lt P(f), P(f)\rt (x) = f(x)^*f(x) \le  \sum_{x = s(\gm)} f(\gm)^*f(\gm) 
= \lt f, f \rt (x)$$  for each $x \in \Gm^0$;
hence, $\|P(f)\|^2 \le \|f\|^2$.  Thus, $P$ extends to a projection  
$q \in \Ll(L^2(E))$; note that $q$ is a projection onto a Hilbert 
submodule isomorphic to $C_0(E^0)$.  If $f \in \cc$ is regarded as 
a left multiplication operator then  
$\|P(f)\| = \|qfq\| \le \|f\|$.  Hence, $P$ extends uniquely to a
linear map, also denoted $P$, from $\cse$  to $C_0(E^0)$ which
restricts to the identity on  $C_0(E^0)$  and $P$ is a projection 
of norm one. \qef

{\bf  3.7 Cor:\ }   For  all $f \in \cc$, $\|f\|_2 \le \|f\|$; thus, 
the inclusion, $\cc \subset L^2(E)$,  extends to a continuous 
map, $$\io: \cse \ra L^2(E).$$ 
\pr  For  all $f \in \cc$, $(\|f\|_2)^2 = \|P(f^*f)\| \le \|f^*f\| =
\|f\|^2$.  \qef

{\bf 3.8  Def:\ }  An open subset $U \subset \Gm$  is said to be an
{\em open $\Gm$-set} if the restrictions, $r|_U$  and $s|_U$, are
one-to-one.  An element $f \in \cc$ is said to be a {\em normalizer} if 
 supp\,$f$  is contained in an open  $\Gm$-set; let $\N(E)$  
denote the collection of all normalizers.
Note:  Since $\cc = {\rm span}\,\N(E)$, $\N(E)$ is total in $\cse$ 
and $L^2(E)$.

{\bf 3.9 Fact:\ }  If $g \in C_c(E^0)$ and $f \in \N(E)$, then  
$f^*gf \in C_c(E^0)$.\\[.5ex]
\pr  For each $x \in \Gm^0$ there is at most one $\gm \in \Gm$  so 
that $s(\gm) = x$ and $f(\gm) \ne 0$.  For $\bt \in \Gm$ with  $s(\bt) = x$,  
if $\bt = x$ and there is such a $\gm$, one has:
$$f^*gf(\bt) =  \sum_{\bt =\gm_1\gm_2\gm_3 } f^*(\gm_1)g(\gm_2)f(\gm_3)
= f(\gm)^*g(r(\gm))f(\gm);$$
if $\bt \ne x$  or there is no such $\gm$, each term in the above sum
is zero and one has $f^*gf(\bt) = 0$.  Thus, $f^*gf \in C_c(E^0)$. \qef

{\bf 3.10  Prop:\ } $P : \cse \ra C_c(E^0)$ is faithful; thus,  
$\io : \cse \ra  L^2(E)$ is injective.\\[.5ex]
\pr  We will show that $P(a^*a) \ge 0$ for every $a \in \cse$, 
$a \ge 0$.  
For all $b \in \cse$  and $f \in \N(E)$ one has:
$$f^*P(b)f = P(f^*bf) = \lt bf, f \rt.$$
It suffices to check this for $b \in \cc$  where it follows by 
a calculation similar to that in the above fact.  
Since $\N(E)$ is total in $L^2(E)$, there is $f \in \N(E)$ so 
that $af \ne 0$; by the above with $b = a^*a$:
$$f^*P(a^*a)f = \lt a^*af, f \rt = \lt af, af \rt \ne 0.$$
Hence $P(a^*a) \ne 0$ and $P$ is faithful. \qef

{\bf 3.11 Fact:\ } The norm on $\cc$ is the unique C*-norm extending 
the supremum norm on $C_c(E^0)$ for which $P$  extends to the 
completion as a faithful conditional expectation.\\[.5ex]
\pr  Let $A$  denote the completion of $\cc$  in such a norm. $P$  
extends to a conditional expectation so left multiplication in 
$\cc$  extends to a continuous *-homomorphism, $A \ra \Ll(L^2(E))$; 
since $P$ is faithful the map must be injective.  \qef

{\bf 3.12\  } Let $\Om$  be an open subgroupoid of $\Gm$; denote the 
inclusion map $j : \Om \ra \Gm$ and put $D = j^*(E)$.

{\bf Prop:\ }  The inclusion $C_c(D) \subset \cc$  is isometric and 
thus extends to an inclusion $\cs(D) \subset \cse$.\\[.5ex]
\pr  That the inclusion extends to a *-homomorphism is immediate.  
The norm on $C_c(D)$ may be characterized as the C*-norm which extends
the supremum norm on $C_c(D^0)$ for which $P$ is faithful on the 
completion.  Since $P$ commutes with the inclusion $C_c(D) \subset
\cc$  and $P$ is faithful on the closure of $C_c(D)$ in $\cse$, the 
inclusion is isometric. \qef

{\bf 3.13\ } Observe that $\|f\|_{\infty} \le \|f\|_2 \le \|f\|$ 
for every  $f \in \cc$.

{\bf \S 4  Morita equivalence}

We continue to restrict attention to principal r-discrete groupoids.  
Let $E$ be a saturated Fell bundle over a groupoid, $\Gm$, and $V$ 
be the associated Hilbert $E^0$-module bundle over $\Gm^0$ (see 3.3), 
we show below that there is an action $\sg$ of $\Gm$ on the C*-bundle
$\K(V)$ so that $\cs(\Gm \times_\sg \K(V))$ and $\cse$ are strong 
Morita equivalent.

{\bf 4.1 Def:\ }  A groupoid morphism, $\ph: \Gm \ra \Dl$, is said to 
be full if for every $x \in \Gm^0$  there is $\gm \in \Gm$  with 
$\ph(\gm) \not\in \Dl^0$ such that $s(\gm) = x$.  For $i=0,1$ set  
$\Gm_i = \ph^{-1}(i)$ and note that $\Gm_i$ is an open subgroupoid 
of $\Gm$; note further that $\Gm_0$  and $\Gm_1$ are equivalent 
(see \cite{MRW}).  Let $j_i : \Gm_i \ra \Gm$ denote the embeddings; 
if $E$ is a Fell bundle over $\Gm$, put $E_i = j_i^*(E)$.

{\bf 4.2 Th:\ }  Let $\Gm$ be a groupoid, $\ph: \Gm \ra \Dl$ be a full
groupoid morphism, and $E$ a saturated Fell bundle over $\Gm$, then, 
with notation as above, $\cs(E_0)$  and $\cs(E_1)$ are strong Morita 
equivalent (cf.\ \cite[Cor. 5.4]{Re2}).\\[.5ex]
\pr  By Prop. 3.12 we may regard $\cs(E_0)$  and $\cs(E_1)$  as
subalgebras of $\cse$.  We show below that $\cs(E_0)$  and $\cs(E_1)$ 
are complementary full corners in $\cse$; it will then follow 
(see \cite[Th. 1.1]{BGR}) that they are strong Morita equivalent.  
It is clear that these subalgebras are complementary corners (as 
in 3.4 the characteristic functions on the unit spaces of $\Gm_0$  
and $\Gm_1$  may be identified with complementary projections in  
$M(\cse)$) and so it remains to show that they are full.  By symmetry 
we need only show that $\cs(E_0)$  is contained in the ideal generated 
by  $\cs(E_1)$.  It suffices to show that $C_c((E_0)^0)$ is 
contained in this ideal.  Since $\ph$  is full, for every 
$x \in (\Gm_0)^0$  there is  $\gm \in \Gm$ with $\ph(\gm) = \df$ 
 and $s(\gm) = x$.  Choose an open $\Gm$-set $U$ containing  $\gm$ 
so that $U \subset \ph^{-1}(\df)$; every $g \in \cc$ with 
supp\,$g \subset \ph^{-1}(\df)$ is in the ideal generated by $\cs(E_1)$,
since  supp\,$gg^* \subset \ph^{-1}(\df\df^*) = \Gm_1$.  Since $E$ is 
saturated we may regard the restriction of $E$ to $U$ as a full
Hilbert module bundle over the restriction of $(E_0)^0$ to $s(U)$.  
Hence, given $f \in C_c((E_0)^0)$  with supp\,$f \in s(U)$  and 
$\ep > 0$, there are $g_k, h_k \in \cc$ with  supp\,$g_k$, supp\,$h_k 
\subset U$ for  $k = 1,\ldots, n$ such that:
$$\| f - \sum_{1 \le k \le n} g_k^*h_k \| < \ep.$$
Hence $f$ is in the ideal generated by $\cs(E_1)$.  Since each
element in $C_c((E_0)^0)$ can be written as the finite sum of such 
elements, $C_c((E_0)^0)$ is contained in the ideal as desired. \qef

{\bf 4.3\ }  Let $E$ be a saturated Fell bundle over $\Gm$ and $V$ be 
the associated Hilbert $E_0$-module bundle over $\Gm_0$; we construct 
another Fell bundle $F$ using $V$; for  $\gm \in \Gm$ set:
$$F_\gm = V_{r(\gm)} \ot_ {E_{r(\gm)}} E_\gm 
\ot_ {E_{s(\gm)}} V_{s(\gm)}^*;$$
note that $F_x = V_x \ot_ {E_x} V_x^* = \K(V_x)$.  Involution is
defined in the obvious way:
$$ u \ot e \ot v^* \mapsto v \ot e^* \ot u^*;$$
given  $(\al, \bt) \in \Gm^2$; if $ t \ot d \ot u^* \in F_\al$ and 
$ v \ot e \ot w^* \in F_\bt$, define multiplication by the formula:
$$(t \ot d \ot u^*)(v \ot e \ot w^*) =
t \ot d\lt u, v \rt e \ot w^* \in  F_{\al\bt}.$$
The verification of the Fell bundle properties is straightforward 
(associativity follows from the associativity of the tensor product 
of equivalence bimodules - see [Ri1, Prop. 6.21]).  Since $E$ is 
saturated and $V$ is full, $F$ is saturated.

{\bf 4.4\ } The following proposition is analogous to \cite[Th.\ 5.3]{Yn} 
in which the existence of an action, which is then defined to be the 
dual action of a given coaction (see \cite[Def.\ 5.4]{Yn}), is 
established.

{\bf Prop:\ }  Let $F$ be as above.  There is an action,  
$\sg : \Gm * \K(V) \ra \K(V)$, so that  
$F \cong \Gm \times_\sg \K(V)$.\\[.5ex]
\pr  First we identify $F$ with the pull back bundle, $\Gm * \K(V)$, 
by means of the following (see 3.3):
$$F_\gm = (V_{r(\gm)} \ot_ {E_{r(\gm)}} E_\gm)
\ot_ {E_{s(\gm)}} V_{s(\gm)}^*\cong 
V_{s(\gm)} \ot_ {E_{s(\gm)}} V_{s(\gm)}^*\cong \K(V_{s(\gm)});$$
the action $\sg_\gm : \K(V_{s(\gm)}) \ra \K(V_{s(\gm)})$ is likewise defined 
by the isomorphism (see also 3.5iv):
$$\K(V_{s(\gm)}) \cong F_\gm = V_{r(\gm)} \ot_ {E_{r(\gm)}} 
(E_\gm \ot_ {E_{s(\gm)}} V_{s(\gm)}^*) \cong
V_{r(\gm)} \ot_ {E_{r(\gm)}} V_{r(\gm)}^*\cong \K(V_{r(\gm)}).$$
On elementary tensors of the form $ac \ot b^* \in \K(V_{s(\gm)})$, 
where $a \in E_\al$, $b \in E_\bt$, $c \in E_\gm$, and $s(\al) = r(\gm)$, 
 $s(\bt) = s(\gm)$, (note that $ac \in E_{\al\gm } \subset V_{s(\gm)}$
and $b \in E_\bt \subset V_{s(\gm)}$ and that elements of the form, 
$ac \ot b^*$ span a dense subset of $\K(V_{s(\gm)})$), $\sg_\gm$ is 
given by:
$$\sg_\gm(ac \ot b^*) = a \ot cb^*.$$
It is a routine matter to verify that $\sg$ defines an action of  
$\Gm$ on $\K(V)$ (for example,  $\sg_x = {\rm id}_{\K(V_x)}$  for   
$x \in \Gm^0$ follows from the fact that 
$\K(V_x) = V_x \ot_{E_x} V_x^*$  is a balanced tensor product) 
and that $F \cong \Gm \times_\sg \K(V)$. \qef

{\bf 4.5 Cor:\ }  With notation as above $\cse$ and 
$\cs(\Gm \times_\sg \K(V))$  are strong Morita equivalent.\\[.5ex]
\pr To apply Th.\ 4.2, we construct a Fell bundle $D$ over 
$\Gm \times \Dl$ which restricts to $E$ on $\Gm \times 0$ and to $F$ 
on $\Gm \times 1$ (by Prop.\ 4.4,  $F \cong \Gm \times_\sg \K(V)$).  
For $\gm \in \Gm$ define:
 $D_{(\gm, 0)} = E_\gm$, 
$D_{(\gm, 1)} = F_\gm$, 
$D_{(\gm, \df)} = V_{r(\gm)} \ot_ {E_{r(\gm)}} E_\gm$, and
$D_{(\gm, \df^*)} = E_\gm \ot_ {E_{s(\gm)}} V_{s(\gm)}^*$.
One defines multiplication and involution in the natural way.  
For example, if $(\al, \bt) \in \Gm^2$  the map:
$$ D_{(\al, \df^*)} \times D_{(\bt, \df)} \ra 
D_{(\al\bt, 0)} = E_{\al\bt}$$
is given by the formula:
$$(d \ot u^*)(v \ot e) =  d\lt u, v \rt e;$$
or the map:
$$ D_{(\al, \df)} \times D_{(\bt, \df^*)} \ra D_{(\al\bt, 1)} = F_{\al\bt} 
= V_{r(\al)} \ot_ {E_{r(\al)}} E_{\al\bt} \ot_ {E_{s(\bt)}} V_{s(\bt)}^*,$$
is given by the formula:
$$(v \ot e)(d \ot u^*) = v \ot ed \ot u^*.$$
Note that $D$ is saturated.  Thus, the theorem applies (the map,  
$\Gm \times \Dl \ra \Dl$, is given by projection onto the second
factor). \qef

{\bf 4.6 Remark:\ } If the groupoid is a (discrete) group this result 
may be obtained by combining \cite[Cor. 2.7]{Q2} and \cite[Th. 8]{Kt}; 
I wish to thank Quigg for bringing this to my attention.

Department of Mathematics,
University of Nevada,
Reno, NV 89557

{\tt alex@math.unr.edu}

\end{document}